\documentclass{article}

\usepackage{arxiv}

\usepackage[nolist]{acronym}
\usepackage{color}
\usepackage{soul}
\usepackage[colorlinks=false,citecolor=red,linkcolor=black]{hyperref}

\usepackage{amsmath}
\usepackage{amsfonts}
\usepackage{amssymb}
\usepackage{amsbsy}

\usepackage[utf8]{inputenc} 
\usepackage[T1]{fontenc}    
\usepackage{url}            
\usepackage{booktabs}       
\usepackage{amsfonts}       
\usepackage{nicefrac}       
\usepackage{microtype}      
\usepackage{cleveref}       
\usepackage{lipsum}         
\usepackage{graphicx}
\usepackage{natbib}
\usepackage{doi}

\usepackage{adjustbox}

\DeclareMathAlphabet\mathbfcal{OMS}{cmsy}{b}{n}

\title{Optimal Scheduling of Variable Speed Pumps using Mixed Integer Linear Programming - Towards an Automated Approach}

\author{
	Tomasz Janus
	\thanks{Visiting Researcher} \\
	Department of Mechanical, Aerospace and Civil Engineering\\
	The University of Manchester\\
	Sackville Street, Manchester, M13 9PL, UK\\
	\texttt{tomasz.janus@manchester.ac.uk} \\
	\And
	Bogumil Ulanicki
	\thanks{Emeritus Professor of Engineering Systems} \\
	School of Engineering and Sustainable Development\\
	De Montfort University\\
	The Gateway, Leicester LE1 9BH, UK\\
	\texttt{bul@dmu.ac.uk} \\
	\And
	Kegong Diao
	\thanks{Senior Lecturer} \\
	School of Engineering and Sustainable Development\\
	De Montfort University\\
	The Gateway, Leicester LE1 9BH, UK\\
	\texttt{kegong.diao@dmu.ac.uk} \\
}

\hypersetup{
pdftitle={Optimal Scheduling of Variable Speed Pumps in Water Distribution Systems using Mixed Integer Linear Programming},
pdfsubject={math.OC},
pdfauthor={Bogumil Ulanicki, Tomasz Janus},
pdfkeywords={milp, energy optimization, water distribution networks},
}

\begin{document}
\maketitle

\begin{abstract}
	This article describes the methodology for formulating and solving optimal pump scheduling problems with \acp{VSP} as \acp{MILP} using piece-linear approximations of the network components.
	The \ac{WDN} is simulated with an initial pump schedule for a defined time horizon, e.g. 24 hours, using a nonlinear algebraic solver.
	Next, the network element equations including \acp{VSP} are approximated with linear and piece-linear functions around chosen operating point(s).
	Finally, a fully parameterized \ac{MILP} is formulated in which the objective is the total pumping cost.
	The method was programmed in MATLAB/OCTAVE and Python and is publicly available on GitHub \footnote{\href{https://github.com/tomjanus/milp-scheduling}{https://github.com/tomjanus/milp-scheduling}}.
	The method was used to solve a pump scheduling problem on a a simple two variable speed pump single-tank network that allows the reader to easily understand how the methodology works and how it is applied in practice.
    The obtained results showed that the formulation is robust and the optimizer is able to return global optimal result in a reliable manner for a range of operating points.
    The work summarized here is a prototype of a framework that is being implemented in a Python package for automated solution of optimal pump scheduling problems on EPANET networks using mixed integer linear programming.
\end{abstract}


\section{Introduction}
\label{sec:introduction}

Pump scheduling plays a crucial role in optimizing \ac{WDN} operation and has been an important subject of research over the last decades.
The purpose of pump scheduling is to find the sequence of pump ON/OFF statuses and pump speeds, in case of \acfp{VSP}, that minimize the pumping cost by shifting pumping to time periods with lower electricity tariffs and by routing the flow through the network such that the energy is used most efficiently.
Optimal pump schedules depend on the electricity tariff profile, the demand profile, and the hydraulic characteristic of the network. 
Pump scheduling is an important problem in the operation of \acp{WDN} because up to 70\% of total operating costs are attributed to electricity consumption for pumping.
Past research and industrial case studies have shown that optimization of pump schedules can lead to up to 10-20\% reduction in pumping costs, i.e. up to 7\%-14\% of total operating costs of \acp{WDN}.

Nevertheless, finding optimal pump schedules is a difficult dynamic mixed-integer nonlinear nonconvex optimization problem \citep{Bonvin2017}.
Nonconvexity arises from the nonlinearities imposed by the network equations that act as constraints.
Integer variables used for selecting statuses of pumps and other active network components such as valves, as well as, in case of piece-wise approximations and relaxations, selection of active domains, make the problem combinatorial.
Dynamics arise in storage tanks and require the optimization problem to consider variables and constraints from all time-steps, which in turn substantially increases the problem size.
Consequently, \acp{PSP} are inherently difficult to solve and, depending on the size of the model and on the optimization method, may get stuck at local minima, fail to find a global optimum within the allocated time or at all, or fail to find a feasible solution altogether.

Optimization methods, including those applied to pump scheduling, can be broadly classified into two categories: (a) deterministic mathematical programming, such as mixed integer nonlinear programming, mixed integer linear programming, etc. and (b) stochastic evolutionary searches, such as e.g. \acp{GA}, \ac{PSO}, and \ac{ACO}. 
The former use full or partial information about the model to guide the optimization process whilst the latter explore the objective space without reliance on the information about the model.
Consequently, mathematical programming tends to be significantly faster at the cost of being more difficult to formulate, while evolutionary algorithms tend to be slow but relatively easy to set up and able to work in conjunction with an input-output model of any mathematical form.
Slower convergence speeds however, are leveraged by some \acp{GA}' massive parallelization capabilities, as demonstrated in \citet{Reed2014}.
Additionally, many \acp{GA} support multiple objectives while multi-objective extensions within mathematical programming frameworks are far less common.
Last but not least, \acp{EA}, including \acp{GA}, are very good at exploring large decision spaces.
Unfortunately, they cannot guarantee that the global optimal solution has been found nor provide bounds on global optimality \citep{Menke2016} which is in contrast to (convex) mathematical programming which can provide such guarantees.

%

The literature on the subject of pump scheduling is extensive and voluminous.
For a more in-depth study the readers are referred to the most recent review papers on the topic by \citet{Jetmarova2017} and \citet{Wu2018optimal}.
Here, we shall only mention a handful of aspects of pump scheduling that pertain to mixed integer programming which has recently become more popular, most likely in response to the recent improvements in speed, reliability, scalability and stability of numerical solvers, such as CPLEX \citep{cplex2009v12}, GUROBI \citep{gurobi}, MOSEK \citep{mosek}, SCIP \citep{BestuzhevaEtal2021OO}

To increase the speed and robustness of finding globally optimal solution, the originally non-convex \ac{MINLP} formulation of the pump scheduling problem needs to be numerically simplified.
The literature distinguishes between (a) model simplification/reduction that is applied to the model before optimization problem formulation and (b) simplification of the optimization problem using various mathematical techniques from the \ac{OR} field.
Model network simplifications were described in \citet{Anderson1995, Deuerlein2008, Alzamora2014} and in the context of real-time pump scheduling, by \citet{Shamir2008}.
The most popular methods for problem simplification are via (a) various relaxations and approximations of constraints, objective(s) and the type of decision variables in order to turn the original problem it into a convex nonlinear or linear problem, (b) decomposition into several easier to solve problems, and (c) relaxation of the optimality criterion see e.g. \cite{Gleixner2012}.

Non-convexity was addressed in \citet{fooladivanda2018, singh2019optimal} who turned the non-convex \ac{MINLP} problem into a convex \ac{MINLP} problem via different \ac{SOC} relaxations.
Additionally, \citet{fooladivanda2018} incorporated \acp{VSP} and \acp{PRV}.
\citet{Bonvin2017} addressed the non-convexity in a less formal and more heuristic way that required restricted formulations supporting only the network topologies without loops.

Problem decomposition into short-term and long-term optimization was investigated in \citet{PulidoCalvo2011}.
Similarly, \citet{Ulanicki2007} proposed time decomposition via solution of a relaxed continuous problem to find optimum reservoir trajectories, followed by a solution of a mixed-integer pump scheduling problem that tracks those trajectories.
Lagrangian decomposition and Benders decomposition were successfully applied to \ac{MINLP} problems by \citet{Ghaddar2015} and \citet{Sawaya2015}, respectively. 
Alternatively, the problem could also be decomposed spatially by dividing the network into smaller isolated subsystems.


Various mixed-integer problem reductions via approximations and relaxations were recently performed by a number of authors.
\citet{vieira2020} used an iterative approach with a feedback loop from EPANET simulator to iteratively limit the error from \ac{MILP} with component relaxations.
In a conceptually similar way, \citet{Liu2020} created a \ac{MILP} formulation with component relaxations which are tightened by the solutions of a series of EPANET simulations.
The authors claimed that their formulation with pipe characteristic relaxations over-performs one with piece-linear approximations due to reduction of binary variables.
\citet{Salomons2020} tested different levels of reduction of binary variables in order to align computational times of \ac{MILP} pump scheduling schemes for real-time pump scheduling applications.
\citet{Bonvin2021} developed relaxation of non-convex constraints of the original non-convex \ac{MINLP} problem using Polyhedral Outer Approximations (OA) and solved the relaxed convex problem with branch and bound method for convex \acp{MINLP}.
Their method supports \acp{VSP}.
Most recently, \citet{tasseff2022} developed tight polyhedral relaxations of the original \ac{MINLP}, derived novel cuts using duality theory, added novel optimization-based bound tightening and cut generation procedures and implemented their method in an open-source Julia package \citep{taseff2019}
The authors addressed two main deficiencies of current state-of-the-art \ac{MILP} solvers: the slow improvement of dual bounds and the difficulty in generating feasible primal solutions.
The authors considered \acp{FSP} only.

Although many advancements in pump scheduling using mixed integer programming have been introduced recently, the literature on the subject is still fragmented with papers addressing some crucial aspects of the methodology whilst omitting others.
Additionally, some treatment of crucial network components such as \acp{VSP} is under-represented.
Meanwhile, the findings reported by many authors suggest that current state-of-the-art \ac{MILP} solvers are able to find optimal pump schedules for networks with around $100$ nodes \cite{Liu2020}, i.e. of sufficient complexity to make them practical
This suggests that an automated method for solving pump scheduling problems with mixed integer programming could be of practical value to the community.
In this paper we communicate our initial findings in prototyping a method for automatic network conversion into a \ac{MILP} problem and subsequent solution using one of the available solvers.
We provide a complete mathematical description of the \ac{MILP} formulation of a network with \acp{VSP}, which is novel.
The validity of the approach and its reliability and robustness is tested on a small network with two parallel \acp{VSP} and one tank.
This work is a part of a larger project on automating pump scheduling with mixed integer linear programming that is being developed in the \verb|dev-python| branch of the GitHub repository of MILOPS-WDN - the Mixed Integer Linear Optimal Pump Scheduler \citep{milops-wdn2023}.
The prototype source code used in this study is contained in the \verb|main| branch of the same repository.

\section{Methodology}
\label{sec:methodology}
\begin{figure}[!b]
	\centering
	\includegraphics[width=0.95\linewidth]{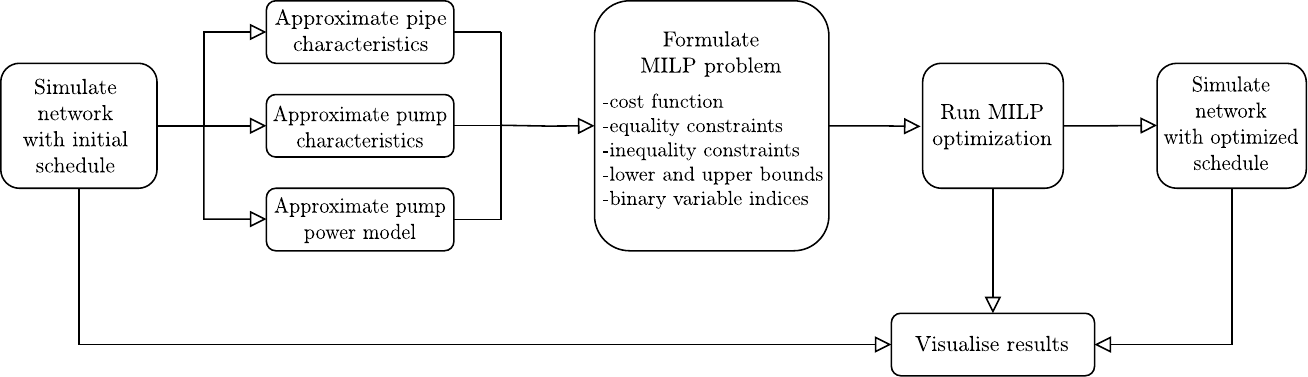}
	\caption{Block diagram visualising the methodology for formulating and solving pump scheduling problems using mixed integer linear programming approach with linear network element approximations.}
	\label{fig:methodology_diagram}
\end{figure}
The procedure described in this paper follows the methodology shown in Fig.~\ref{fig:methodology_diagram}.
First, the network is simulated using initial schedules in order to obtain approximate operating points that are used in the subsequent approximations of the model components and the pump power consumption characteristic.
The nonlinear equations of pipe and pump characteristics are approximated using linear and piecewise linear approximations so that the network model can be put into a linear form required by \ac{MILP}.
Other nonlinear components such as valves, including \acp{CV} and \acp{PRV}, as well as other nonlinearities such as leakage could additionally be included in the model formulation.
For simplicity, these are not considered in this study, but the methodology does not prohibit their inclusion.
Subsequently, the \ac{MILP} optimization problem of the standard form shown in Eq.~\ref{eq:milp_formulation} is solved for the desired time horizon of typically $24$ hours.
Finally, the optimal schedule of pump ON/OFF statuses and pump speeds is input back into the simulator and the final simulation result is compared against the result of the initial simulation.

\begin{equation}
\begin{aligned}
\min_{\mathbf{x}} \quad & \mathbf{c}^T\,\mathbf{x}\\
\textrm{s.t.} \quad & \mathbf{A}_{ineq}\,\mathbf{x} \le \mathbf{b}_{ineq}\\
  & \mathbf{A}_{eq}\,\mathbf{x} = \mathbf{b}_{eq}\\
  & \mathbf{l} \le \mathbf{x} \le \mathbf{u}\\
  &x_i \in \mathbb{Z}, \forall_i \in \varUpsilon   \\
\end{aligned}
\label{eq:milp_formulation}
\end{equation}
The objective function $\mathbf{c}^T\,\mathbf{x}$describes the total pumping cost over the optimization time-horizon.
The network equations are described with inequality and equality constraints using with two tuples: $\left(\mathbf{A}_{ineq}, \mathbf{b}_{ineq}\right)$ and $\left(\mathbf{A}_{eq}, \mathbf{b}_{eq}\right)$.
The lower bounds $\mathbf{l}$ and the upper bounds $\mathbf{u}$ enforce physical limits on the decision variable vector $\mathbf{x}$, such as e.g. minimum and maximum tank levels, maximum pump flows, etc. $\varUpsilon$ is a nonempty subset of the set $\{1\ldots n\}$ that specifies the indices of the integer variables, where $n = |\mathbf{x}|$. Integer variables are used for the selection of pumps and active segments in piece-linear approximations.
The decision variable vector $\mathbf{x}$ includes all state variables of the network model in all time moments plus the auxiliary (artificial) variables.
The auxiliary variables are introduced to represent the variables bound to the domains of piece-linear segments obtained via piece-linear linear approximations.
By convention, an auxiliary (continuous) variable for $x$ has a symbol $xx$ with an index representing the index of the subdomain, e.g. $xx_i$ represents the value of $x$ if $x$ lies inside segment $i$.
By definition $x = \sum_{i=1}^m xx_i$ where $m$ is the number of segments.
Auxiliary variables are accompanied by binary selection integer variables.
By definition, these are denoted by capital letters, e.g. $XX$, and satisfy equation $\sum_{i=1}^m XX_i = 1$, which means that only one component can be active at a time.

\subsection{Network equations}
\label{subsec:network_equations}
Water network is described with standard network equations describing: (1) headlosses in pipes, (2) flow continuity equations in nodes, incl. tanks, (3) pump characteristics describing pumping head and power consumption vs. flow, and (4) head vs. volume relationships in tanks. More information about modelling of \acp{WDN} can be found in \citet{strafaci2007advanced}. The equations are briefly listed below as their understanding is necessary to follow the model approximation steps.

\subsubsection{Pipe headloss equations}
\label{subsubsec:pipe_headlosses}
The nodal pressures in the network are modelled as follows:
\begin{equation}
 \mathbf{R} \left|\mathbf{q}(k)\right| \mathbf{q}(k) + \mathbf{\Lambda}_c^T \mathbf{h}_c(k) + \mathbf{\Lambda}_f^T\, \mathbf{h}_f(k) = \mathbf{0}
\end{equation}
where $\mathbf{\Lambda}_c$ and $\mathbf{\Lambda}_f$ are node-element incidence matrices for the connection nodes and the fixed nodes, respectively. $\mathbf{R}$ are pipe resistances, $\mathbf{q}$ is the vector of element flows, $\mathbf{h}_c$ are the calculated heads and $\mathbf{h}_f$ are the fixed heads, e.g. heads in tanks and reservoirs.
The above equation can be split into multiple equations, each representing a single pipe $j$ for $j \in \langle1,\ldots,|\mathbb{P}|\rangle$ where $\mathbb{P}$ denotes the set of pipes.
\begin{equation}
 \underbrace{R_j \left|q_j(k)\right| q_j(k)}_\textrm{pipe characteristic} + \mathbf{\Lambda}_{c,j}^T \mathbf{h}_c(k) + \mathbf{\Lambda}_{f,j}^T\, \mathbf{h}_f(k) = 0
\end{equation}
where $\mathbf{\Lambda}_{c,j}^T \mathbf{h}_c(k) = h_{d,j}(k)$ and $\mathbf{\Lambda}_{f,j}^T\, \mathbf{h}_f(k) = h_{o,j}(k)$, i.e. the downstream and the upstream head of pipe $j$, respectively.

\subsubsection{Mass balance in nodes}
\label{subsubsec:mass_balance_nodes}
Mass balance in nodes for each time step $k \in {1 \ldots K}$ is calculated as
\begin{equation}
 \label{eq:mass_balance_nodes}
 \mathbf{\Lambda}_c\,\mathbf{q}(k) - \mathbf{d}(k) = \mathbf{0}
\end{equation}
where $\mathbf{d}$ is the matrix of nodal demands.

\subsubsection{Pump power consumption}
\label{subsubsec:pump_power_consumption}
Pump power consumption is modelled with the relationship described in \citet{Ulanicki2008} representing power demand of a group of $n$ identical pumps, each operating at speed $s$.
\begin{equation}
 \label{eq:pump_power_consumption_1}
 P(q,n,s) = n\,s^3\,P\left(\frac{q}{n\,s}\right)
\end{equation}
where
\begin{equation}
 \label{eq:pump_power_consumption_2}
 P\left(\frac{q}{n\,s}\right) = a_3\left(\frac{q}{n\,s}\right)^3 + a_2 \left(\frac{q}{n\,s}\right)^2 + a_1 \left(\frac{q}{n\,s}\right) + a_0
\end{equation}
The coefficients $a_3$, $a_2$, $a_1$, $a_0$ are unique for each individual pump model.
Smaller individual \acp{VSP} could alternatively be modelled with the scaling of \citet{Sarbu_Borza_1998} as advised in \cite{Simpson2013}, although this decision is left to the user.

\subsubsection{Pump hydraulics}
\label{subsubsec:pump_hydraulics}
Pump hydraulics are formulated with the pump characteristic model $H=H(q, n, s)$ from \citet{Ulanicki2008} that describes the relationship between the head gain $H$ and the pump group flow $q$ for a group of $n$ identical pumps operating at speed $s$.
\begin{equation}
 \label{eq:pump_head_n_s}
 \frac{H}{n^2\,s^2} = A \, \left( \frac{q}{n\,s} \right)^2 + B \left( \frac{q}{n\,s} \right) + C
\end{equation}
which translates into
\begin{equation}
 \label{eq:pump_head_n_s_simplified}
 H = A \, q^2 + B \, q \, n \, s + C \, n^2 \, s^2
\end{equation}

Eqs.~\ref{eq:pump_power_consumption_1}, \ref{eq:pump_power_consumption_2}, and \ref{eq:pump_head_n_s_simplified} are implemented in the \verb|dev-battery| branch of the EPANET repository \citep{janus2020}.

\subsubsection{Tank model}
\label{subsubsec:tank_model}
Tanks are modelled with a backward finite-difference equation describing the change in tank head $h_t$ between time steps $k$ and $k-1$ as a function of the net flow $q_t(k)$ in/out of the tank in time steps $k = 1 \ldots K$, where $K$ is the optimization time horizon.
\begin{equation}
\label{eq:tank_mass_balance}
h_t(k)-h_t(k-1) - \frac{1}{A_t} \, q_t(k) = 0  \quad \forall k = 2\ldots K
\end{equation}
with initial condition
\begin{equation}
 \label{eq:tank_mass_balance_ic}
 h_t(1) - h_{t,init} = 0
\end{equation}
For cylindrical tanks, the tank's surface area $A_t = \textrm{const}$.
For other geometries, $A_t = A_t(h_t)$ and needs to be included in the model.
If this relationship is not linear, it needs to be approximated with one of the approximation methods - see below.

\subsection{Linear and piece-linear approximations of the network components}
\label{subsec:linearizations}
When approximating nonlinear functions with linear functions, the choice of the linearization technique is often left to the user.
The possible choices are: (1) tangent line approximation using first two terms of Taylor's expansion, (2) linearization by substitution via introduction of additional variables and transformations, (3) piecewise linearization, (4) convex hull approximation.
The choice of the method should be based on the specific characteristics of the nonlinear function and the context of the problem.
Each method has its own limitations and applicability and, in the context of pump scheduling, will affect the accuracy of the solution and the complexity of the \ac{MILP} formulation.
In the following sections, the choice of linearization techniques was made by the Authors but the methodology is not limited to those choices and the readers are encouraged to try different techniques in order to fine-tune their problem formulations.

\subsubsection{Linear approximation of the pump power consumption model}
\label{subsubsec:linearized_pump_power_consumption}
Pump power consumption is linearized by finding a tangent line to the power consumption model of the group of $n$ parallel pumps given in Eq.~\ref{eq:pump_power_consumption_1} and Eq.~\ref{eq:pump_power_consumption_2} at the linearization point $(q_0, s_0)$.
Since each pump is linearized individually, the linearized equations are derived for $n=1$.
Linearization using a tangent at the nominal point is chosen over other linearization methods due to the fact that power consumption curves tend to be quite flat.
Therefore, linearizing the curve around the nominal operating point should offer sufficient approximation accuracy whilst keeping the complexity at minimum, e.g. compared to piece-wise linearization.
Linearization of Eq.~\ref{eq:pump_power_consumption_1} and Eq.~\ref{eq:pump_power_consumption_2} for $n=1$ around the selected operating point $(q_0, s_0)$ yields
\begin{equation}
 P(q,s) = P(q_0,s_0) + 3 \, a_3 \, q_0^2 \, \delta q \, + 2 \, a_2 \, q_0 \, s_0 \, \delta q + a_2 \, q_0^2 \, \delta s \, + a_1 \, s_0^2 \, \delta q + 2 \, a_1 \, q_0 \, s_0 \, \delta s \, + 3 \, a_0 \, s_0^2 \, \delta s
\end{equation}
After grouping similar terms with $\delta s$ and $\delta q$
\begin{equation}
 P(q,s) = P(q_0,s_0) + \left( 3 \, a_3 \, q_0^2 + 2 \, a_2 \, q_0 \, s_0 + a_1 \, s_0^2\right) \delta q + \left( a_2 \, q_0^2 + 2 \, a_1 \, q_0 \, s_0 + 3 \, a_0 \, s_0^2 \right) \delta s
\end{equation}
where $\delta q = q - q_0$ and $\delta s = s - s_0$. Using short notation $P(q,s) = P$ and $P(q_0, s_0) = P_0$ we obtain the following equation for the difference between the power consumption at point $(q, s)$ and the power consumption at the linearization point $(q_0, s_0)$.
\begin{equation}
 \label{linearized_pump_power_long}
 P' = P - P_0 = \left( 3 \, a_3 \, q_0^2 + 2 \, a_2 \, q_0 \, s_0 + a_1 \, s_0^2\right) \, (q - q_0) + \left( a_2 \, q_0^2 + 2 \, a_1 \, q_0 \, s_0 + 3 \, a_0 \, s_0^2 \right) \, (s - s_0)
\end{equation}
Ultimately, the linearized power consumption of a pump at any time step $k$ is calculated with Eq.~\ref{eq:linearized_pump_power}.
\begin{equation}
 \label{eq:linearized_pump_power}
 P(q(k), s(k)) = P(q_0, s_0) + P'(q(k), s(k))
\end{equation}
Throughout this paper $(q_0, s_0) = (q_n, s_n)$, i.e. the power is linearized around the operating point at the nominal speed $s_n=1.0$ and flow $q_n$ equal to the flow at the nominal speed for which the pump efficiency is at its maximum.
After collecting all terms, Eq.~\ref{linearized_pump_power_long} can be simplified to Eq.~\ref{eq:linearized_pump_power_simplified}. The expressions for coefficients $m^q_j$, $m^s_j$ and $c_j$ can be found in the Appendix in Section~\ref{app:lin_pump_power_expressions}.
\begin{equation}
 \label{eq:linearized_pump_power_simplified}
 P_j(k) = m^q_j \, q_j(k) + m^s_j \,s_j(k) + c_j
\end{equation}
Eq,~\ref{eq:linearized_pump_power_simplified} holds only if the pump status is ON.
This is enforced in \ac{MILP} by expressing the equality as a double sided inequality with the `big U' trick as expressed in Eq.~\ref{eq:linearized_pump_power_inequality}.
$U_{power}$ is a large number in the order of magnitude but larger in value than the largest power consumption calculated by the model and $n_j(k)$ is the status of $j$-th pump at timestep $k$.
\begin{equation}
 \label{eq:linearized_pump_power_inequality}
 (n_j(k)-1) \, U_{power} \le m^q_j \, q_j(k) + m^s_j \,s_j(k) + c_j - P_j(k) \le (1-n_j(k)) \, U_{power}
\end{equation}
If $n_j(k) = 1$, Eq.~\ref{eq:linearized_pump_power_inequality} is reduced to equality as both sides of the inequality are zero.
If $n_j(k) = 0$, $(n_j(k)-1) \, U_{power} = -U_{power}$ and $(1-n_j(k)) \, U_{power} = +U_{power}$. Consequently, the equality in Eq.~\ref{eq:linearized_pump_power_simplified} is not enforced.
Note that Eq.~\ref{eq:linearized_pump_power_simplified} yields non-zero power consumption for $q=0$ and $s=0$ as an unwanted byproduct of linearization.
To ensure that power is null when the pump is OFF, $P_j(k)$ is made to obey the following two-sided inequality that forces $P_j(k) = 0$ for $n_j(k) = 0$.
\begin{equation}
 \label{eq:power_binary_linearization}
 0 \le P_j(k) \le n_j(k) \, U_{power}
\end{equation}


\subsubsection{Piece-linear approximation of the pipe model}
\label{subsubsec:linearized_pipe_model}
\begin{figure}[!h]
	\centering
	\includegraphics[width=0.42\linewidth]{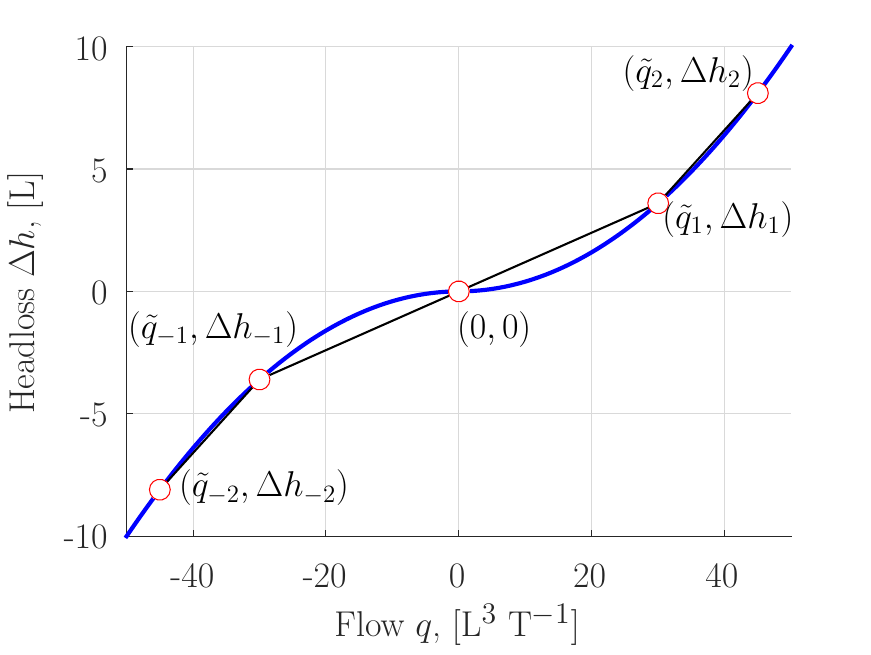}
	\caption{A generic form of a quadratic pipe characteristic model describing headloss $\Delta h$ across the pipe vs flow $q$ and its piece-wise linearization with $3$ linear segments.}
	\label{fig:pipe_linearization}
\end{figure}

Nonlinear pipe characteristic can approximated by $n_s$ piece-linear segments by defining $n_s-1$ breakpoints and exploiting the symmetry of the characteristic around the origin $(0,0)$.
A generic pipe characteristic and its piece-linear form for $n_s = 3$ segments : $[\tilde{q}_{-2}, \tilde{q}_{-1}]$, $[\tilde{q}_{-1}, \tilde{q}_1]$ and $[\tilde{q}_1, \tilde{q}_2]$, is shown in Figure~\ref{fig:pipe_linearization}.
Each segment $i$ of pipe $j$ is described with a linear equation $\Delta h_j=m_{j,i}\,q_j+c_{j,i}$  where
\begin{equation}
 m_{j,i} = \frac{\Delta h_{j,i} - \Delta h_{j,i-1}}{\tilde{q}_{j,i} - \tilde{q}_{j,i-1}}
\end{equation}
and
\begin{equation}
 c_{j,i} = \frac{\Delta h_{j,i-1}\,\tilde{q}_{j,i} - \Delta h_{j,i}\,\tilde{q}_{j,i-1}}{\tilde{q}_{j,i} - \tilde{q}_{j,i-1}}
\end{equation}
The breakpoint $(\tilde{q}_{j,1}, \Delta h_{j,1})$ is fixed for every pipe $j$ and taken from the simulator, e.g. from the model state at 12.00 o'clock.
The breakpoint $(\tilde{q}_{j,2}, \Delta h_{j,2})$ is chosen arbitrarily for every pipe $j$ and should be made large enough to cover the entire range of observed pipe flows.
The points $(\tilde{q}_{j,-1}, \Delta h_{j,-1})$ and $(\tilde{q}_{j,-2}, \Delta h_{j,-2})$ do not need to be calculated but instead, can be derived by exploiting the symmetry of the pipe's characteristic around point $(0,0)$.
In order to represent the three linear segments within one pipe model, two new types of auxiliary variables are required: a continuous variable $ww_{j,i}$ for the flow in pipe $j$ and segment $i$, and a binary variable $BB_{j,i}$ for selecting the active segment of the piece-wise linearized pipe characteristic.
\begin{equation*}
  ww_{j,i} = \left\{
    \begin{array}{ll}
      q_j, & \mbox{if $q_j\in[\tilde{q}_{j,i-1}, \tilde{q}_{j,i}]$}.\\
      0, & \mbox{otherwise}.
    \end{array}
  \right.
\end{equation*}
\begin{equation*}
  BB_{j,i} = \left\{
    \begin{array}{ll}
      1, & \mbox{if $q_j\in[\tilde{q}_{j,i-1}, \tilde{q}_{j,i}]$}.\\
      0, & \mbox{otherwise}.
    \end{array}
  \right.
\end{equation*}
The piece-linear pipe model is represented by the following 4 relationships.
\begin{equation}
 \label{eq:flow_single_pipe}
 BB_{j,i}(k) \, \tilde{q}_{j,i-1} \le ww_{j,i}(k) \le BB_{j,i}(k) \,\tilde{q}_{j,i}
\end{equation}
\begin{equation}
 \label{eq:pipe_segment_flow_equality}
 q_{j}(k)-\sum_{i=1}^{n_{s,pipe}} ww_{j,i}(k)=0
\end{equation}
\begin{equation}
 \label{eq:pipe_segment_selection}
 \sum_{i=1}^{n_{s,pipe}} \, BB_{j,i}(k)=1
\end{equation}
\begin{equation}
 \label{eq:pipe_headloss_linear}
 h_o^j(k)-h_d^j(k)-\sum_{i=1}^{n_{s,pipe}} \left(m_{j,i}\,ww_{j,i}(k) + c_{j,i}\,BB_{j,i}(k)\right)=0
\end{equation}
where
$h_o^j(k)-h_d^j(k) = \Delta h_j(k) = \Lambda_{c,j}^T h_c(k) + \Lambda_{f,j}^T\, h_f(k)$

Eq.~\ref{eq:flow_single_pipe} `binary linearizes' the flow with respect to segment selection. If the binary segment selection variable $BB_{j,i}$ is zero, $ww_{j,i}$ is forced to become zero. Otherwise, $ww_{j,i}$ is bound between $\tilde{q}_{j,i-1}$ and $\tilde{q}_{j,i}$.
Eq.~\ref{eq:pipe_segment_flow_equality} relates the pipe flow variable $q$ to the auxiliary flow variable $ww$. Eq.~\ref{eq:pipe_segment_selection} ascertains that only one segment in the piece-linear equation is active at a time. Finally, Eq.~\ref{eq:pipe_headloss_linear} describes the linear pipe headloss.


\subsubsection{Piece-linear approximation of the pump characteristic - \acp{VSP}}
\label{subsubsec:linearized_pump_characteristic}
\begin{figure}[!h]
	\centering
	\includegraphics[width=0.95\linewidth]{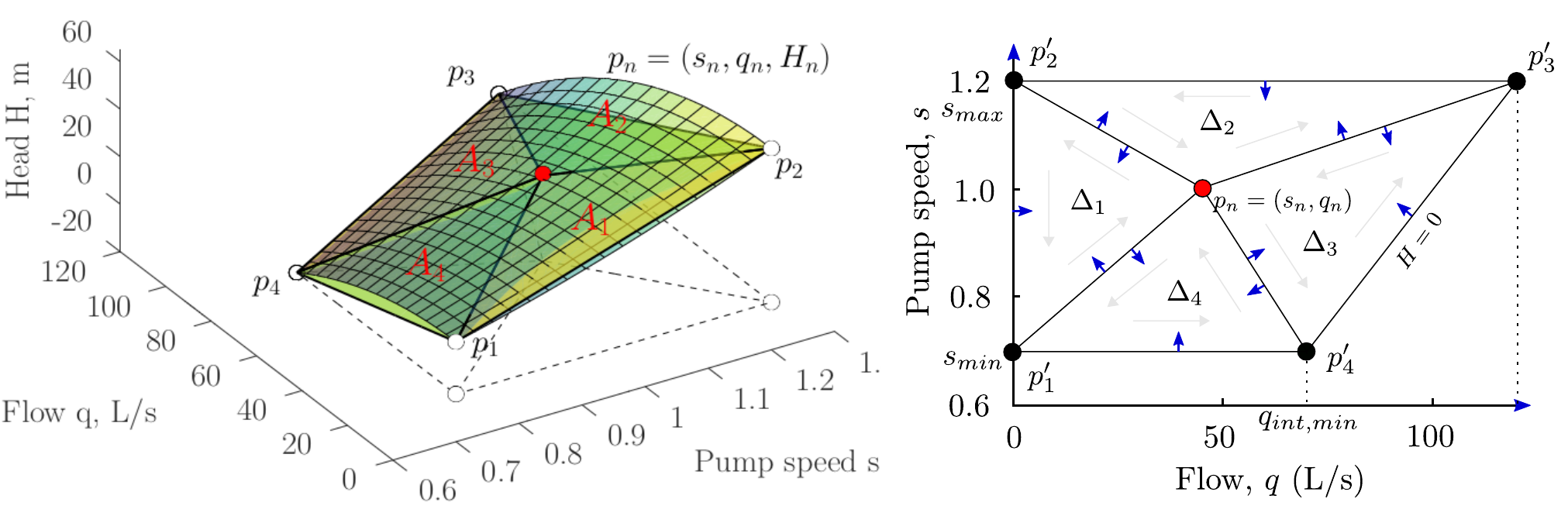}
	\caption{A generic pump characteristic $H = H(q, s)$ with its piece-linear approximation (left) and projection of the piece-linear approximation onto the $q-s$ plane (right).}
	\label{fig:pump_linearization}
\end{figure}

Pump characteristic is approximated with piece-linear surfaces.
A piece linear approximation of a surface can be constructed in different ways with different number of linear segments such as in the case of piece-linear curve linearization, but also with different segment geometries.
In this paper, a piece-linear approximation of the pump characteristic is a top surface of a polyhedron whose sides are defined by vertices $p_1 = (s_{min}, q_{min}, H_1)$, $p_2 = (s_{max}, q_{min}, H_2)$, $p_3 = (s_{max}, q_{int}^{s_{max}}, H_3)$, $p_4 = (s_{min}, q_{int}^{s_{min}}, H_4)$, and the nominal point $p_n = (s_n, q_n, H_n)$ - see Fig.~\ref{fig:pump_linearization}.
The nominal point is derived for the nominal speed $s_n=1.0$ and the maximum efficiency flow $q_n = q^{\eta_{max}}$.
$s_{min}$ and $s_{max}$ are the minimum and the maximum allowed pump speeds, respectively.
$q_{min} = 0$ is the minimum pump flow and $q_{int}^{s_{min}}$ and $q_{int}^{s_{max}}$ are the intercept flows, i.e. the flows for which the pump head $H=0$ at the minimum and at the maximum pump speed, respectively.
Consequently, $H_3 = 0$ and $H_4 = 0$.

The top surface of the polyhedron is defined by the equations of four Euclidean planes $A_i = \{(q, s, H) \, | \, (q, s) \in \Delta_i, H \in \mathbb{R}\}$.
For each plane $A_i$, the point $(q,s)$ must lie within the triangular domain $\Delta_i$ defined by the plane's projection onto the $q-s$ space - see Fig~\ref{fig:pump_linearization} (right).
The equations of each of the four planes $A_i$ are derived from the three vertices of the polyhedron that are contained within it.
The domain for each plane is defined by three inequality constraints derived from the equations of the lines containing the three sides of its projection - see Section~\ref{app:linearized_pump_characteristics_plane_derivation} of the Appendix.
Direction of each inequality is indicated in Fig~\ref{fig:pump_linearization} with a small blue arrow.
The projections of points $p_1$, $p_2$, $p_3$, $p_4$, and $p_n$ onto the $q-s$ space are denoted as $p'_1$, $p'_2$, $p'_3$, $p'_4$, and $p'_n$, respectively.

Introduction of piece-linear pump characteristic requires three auxiliary variables: a binary variable $AA_{j,i}(k)$, and two continuous variables $ss_{j,i}(k)$, $qq_{j,i}(k)$ for each of the four domains, each pump and each time-step.
$AA_{j,i}(k)$ defines whether the current operating point $(s_j(k), q_j(k))$ of pump $j$ lies within the $i$-th segment of the linearized pump characteristic: 1 for YES and 0 for NO.
$ss_{j,i}(k)$ and $qq_{j,i}(k)$ are the speed and the flow of pump $j$ in each domain $i$ of the piece-linear pump characteristic approximation, respectively.
The approximation is defined as follows.
\begin{equation}
 \label{eq:pump_segment_speeds}
 s_j(k) - \sum_{i=1}^{n_{s,pump}} ss_{j,i}(k) = 0
\end{equation}
\begin{equation}
 \label{eq:pump_segment_flows}
 q_j(k) - \sum_{i=1}^{n_{s,pump}} qq_{j,i}(k) = 0
\end{equation}
Eq.~\ref{eq:pump_segment_speeds} and Eq.~\ref{eq:pump_segment_flows} link the auxiliary segment speeds and segment flows to the original pump speeds and pump flows, respectively.
Only one segment is allowed to be active if the pump is ON, i.e. $n_j = 1$
Otherwise, if the pump is OFF, i.e. $n_j = 0$, no segments are allowed to be active.
\begin{equation}
 \label{eq:pump_segment_selection}
 \sum_{i=1}^{n_{s,pump}} AA_{j,i}(k)-n_j(k)=0
\end{equation}
The binary segment selection variable $AA_{j,i}(k)$ is used to `binary linearize' the pump speed and the pump flow with respect to segment selection.
\begin{equation}
  \label{eq:binary_linearization_ss}
  AA_{j,i}(k) \, s_{j,min} \le ss_{j,i}(k) \le AA_{j,i}(k) \, s_{j,max}
\end{equation}
\begin{equation}
  \label{eq:binary_linearization_qq}
  0 \le qq_{j,i}(k) \le AA_{j,i}(k) \, q_{j,max}
\end{equation}
If $AA_{j,i} = 0$, then $ss_{j,i}$ and $qq_{j,i}$ are forced to be zero. Otherwise, $ss_{j,i}$ and $qq_{j,i}$ are bound with box constraints between $s_{j,min}$ and $s_{j,max}$, and $0$ and $q_{j,max} = q_{int,max}$, respectively - see Fig~\ref{fig:pump_linearization} (right).
In summary, pump speeds $s_j(k)$ and pump flows $q_j(k)$ are forced to be equal to the sums of the auxiliary variables $ss_{j,i}(k)$ and $qq_{j,i}(k)$ where, at most one auxiliary segment variable is active at any time-step $k$.
The linearized pump characteristic is represented with the following formula using the `big-U' trick.
\begin{equation}
 \label{eq:linearized_pump_characteristic}
 -U' \le \Delta h_j(k)-\sum_{i=1}^{i=4}\left(dd_{j,i}\,ss_{j,i}(k)+ee_{j,i}\,qq_{j,i}(k)+ff_{j,i}\,AA_{j,i}(k))\right) \le U'
\end{equation}
where $U' = (1-n_j(k)) U_{pump}$ and the coefficients $dd_{j,i}$, $ee_{j,i}$, $ff_{j,i}$ describe the equation of plane $A_i$. Derivation of plane equations is described in Appendix in Section~\ref{app:linearized_pump_characteristics_plane_derivation}.
Each triangular domain $\Delta_i$ of segment $A_i$ is defined with three inequality constraints:
\begin{equation}
\label{eq: pump domain constraints}
\begin{bmatrix}
  m_{qq}^{(1)} & m_{ss}^{(1)} & c^{(1)} \\
  m_{qq}^{(2)} & m_{ss}^{(2)} & c^{(2)} \\
  m_{qq}^{(3)} & m_{ss}^{(3)} & c^{(3)} \\
\end{bmatrix}
\begin{bmatrix}
  qq_{j,i} \\
  ss_{j,i} \\
  1 \\
\end{bmatrix}
\le
\begin{bmatrix}
  0 \\
  0 \\
  0 \\
\end{bmatrix}
\end{equation}

\subsection{\ac{MILP} formulation}
\label{subsec:milp_formulation}
\ac{MILP} formulation of the pump scheduling problem is composed of (1) objective function representing the total pumping cost, (2) a set of equality constraints representing the originally linear and linearized network component equations and auxiliary relationships, (3) a set of inequality constraints representing binary linearized network component equations and additional constraints such as e.g. symmetry-breaking constraints, (4) a set of \ac{LB} and \ac{UB}, aka. box constraints on selected decision variables, and (5) a vector of indices of binary decision variables.

\subsubsection{Objective}

The objective is the total cost of pumping over time horizon $K$, given the energy tariff $T(k)$ and the linearized pumping cost model for each pump $P_j(k)$. $\Delta t$ is the time-step - usually 1h.
\begin{equation}
 \label{eq:total_cost}
 TC = \sum_{k=1}^K \sum_{j\in E_{pump}} P_j (k) \; T(k) \; \Delta t
\end{equation}
%
Note that it is also common for the objective function to include a term penalizing pump switching, e.g. \citep{Lansey1994}.
Since this term includes a sum of absolute (or squared) differences between consecutive pump statuses, it will need to be linearized via introduction of additional variables and constraints \citep{Shanno1971}.
In order to not overcomplicate the current study, inclusion pump switching cost in the objective is left for later.

\subsubsection{Additional constraints}
\paragraph{Symmetry Breaking: }
\label{subsubsec:symmetry_breaking}
Symmetries arise in \ac{MILP} problems when the same feasible solution can be represented in more than one way.
These symmetries can lead to redundant computations as they increase the search space and require the branch \& bound algorithms to explore and compare multiple equivalent branches, slowing down the optimization process.
In pump scheduling, symmetries will arise if parallel pumps within one pumping station have the same characteristic.
These symmetries are removed by introducing an additional set of inequality constraints which enforces the priority of pumps, as described in Eq.~\ref{eq:pump_symmetry_constraint}.
Consequently, the lower priority pumps can be switched ON iff the higher priority pumps are also switched ON, thus preventing the optimizer from needlessly exploring equivalent solutions with different permutations of ON/OFF statuses among equal pump units.
\begin{equation}
 \label{eq:pump_symmetry_constraint}
 -n_{j+1}(k) + n_j(k) \le 0 \quad \forall j \in \{1 \ldots (n_{pumps}-1)\} \quad \textrm{for every pump group with equal pumps}
\end{equation}
Adding this constraint reduces the search space for each pumping station with $n_{pumps}$ from $2^{n_{pumps}}$ to $n_{pumps} + 1$ \citep{Gleixner2012}.

\paragraph{Enforcing tank levels: }
To prevent the optimizer from emptying the reservoirs as it tries to reduce the total pumping cost, the tank level difference between the final time $N$ and the initial time $1$ is bound to a small threshold $\delta_{h_{t}}$ \citep{Menke2016}.
\begin{equation}
 \label{eq:final_tank_level}
	h_{t,j}(N) - h_{t,j}(1) \le \delta_{h_{t,j}} \quad \forall j = 1,\,\ldots\,n_{t}
\end{equation}

The summary of equality and inequality constraints required for the formulation of our pump scheduling problem are listed in Table~\ref{tab:equality_constraints} and Table~\ref{tab:inequality_constraints}, respectively.

\begin{table}[htbp]
 \centering
 \small
 \caption{Equality constraints}
 \begin{tabular}{llcc}
  \toprule
  & Name & Equation(s) & No. of constraints \\
  \midrule
  1 & Mass balance in nodes & \ref{eq:mass_balance_nodes} & $n_n \times K$ \\
  2 & Head-flow relationship in tanks & \ref{eq:tank_mass_balance} + \ref{eq:tank_mass_balance_ic} & $n_t \times K$ \\
  3 & Pipe segment flows & \ref{eq:pipe_segment_flow_equality} & $n_p \times K$ \\
  4 & Pipe segment selection variables & \ref{eq:pipe_segment_selection} & $n_p \times K$ \\
  5 & Linearized pipe headlosses & \ref{eq:pipe_headloss_linear} & $n_p \times K$ \\
  6 & Pump segment speeds & \ref{eq:pump_segment_speeds} & $n_{pump} \times K$ \\
  7 & Pump segment flows & \ref{eq:pump_segment_flows} & $n_{pump} \times K$ \\
  8 & Pump segment selection variables & \ref{eq:pump_segment_selection} & $n_{pump} \times K$\\
  \bottomrule
 \end{tabular}
 \label{tab:equality_constraints}
\end{table}

\begin{table}[htbp]
 \centering
 \small
 \caption{Inequality constraints}
 \begin{tabular}{llcc}
  \toprule
  & Name & Equation & No. of constraints \\
  \midrule
  1 & Pump power binary linearization with respect to pump status &\ref{eq:linearized_pump_power_inequality} & $2 \times n_{pump} \times K$\\
  2 & `Zero power' enforcement for switched OFF pumps & \ref{eq:power_binary_linearization} & $2 \times n_{pump} \times K$ \\
  3 & Binary linearization of pipe flow with respect to pipe segment selection & \ref{eq:flow_single_pipe} & $3 \times n_p \times K$ \\
  4 & Binary linearization of \ac{VSP} speed with respect to pump segment selection & \ref{eq:binary_linearization_ss} & $4 \times n_{pump} \times K$ \\
  5 & Binary linearization of \ac{VSP} flow with respect to pump segment selection & \ref{eq:binary_linearization_qq} & $4 \times n_{pump} \times K$ \\
  6 & Binary linearized \ac{VSP} characteristic & \ref{eq:linearized_pump_characteristic} & $2 \times n_{pump} \times K$ \\
  7 & \ac{VSP} characteristic domain definitions & \ref{eq: pump domain constraints} & $12 \times n_{pump} \times K$ \\
  8 & Symmetry breaking in pump groups with equal pumps & \ref{eq:pump_symmetry_constraint} & $n_{groups} \, \left( n_{eq. pumps}^{group} -1 \right) \times K$ \\
  9 & Enforcing final tank level & \ref{eq:final_tank_level} & $n_{t}$ \\
  \bottomrule
 \end{tabular}
 \label{tab:inequality_constraints}
\end{table}

\subsubsection{Lower and upper bounds on decision variables}
\label{subsubsec:lower_upper_bounds}
Most of the decision variables in vector $\mathbf{x}$ are rather tightly constrained by the inequality and equality constraints.
The exceptions are: (a) heads in tanks, which need to be additionally constrained with box-constraints such that the levels do not violate restrictions imposed by the tanks' minimum ($h^j_{t,min}$) and maximum ($h^j_{t,max}$) levels and (b) integer variables that we restrict to take only binary values.
Tank level constraints for $j \in \{1,\ldots, n_t\}$ tanks, are listed below
\begin{equation}
h^j_{t,min}(k) \le h_t^j(k) \le h^j_{t,max}(k)
\end{equation}
The binary variable constraints can be expressed as follows
\begin{equation}
0 \le x_i \le 1 \quad \forall_i \in \varUpsilon
\end{equation}
where $\varUpsilon$ is the set of indices of integer decision variables - see Eq.~\ref{eq:milp_formulation}


\section{Case study}
\label{sec:case_study}
\begin{figure}[!h]
	\centering
	\includegraphics[width=0.55\linewidth]{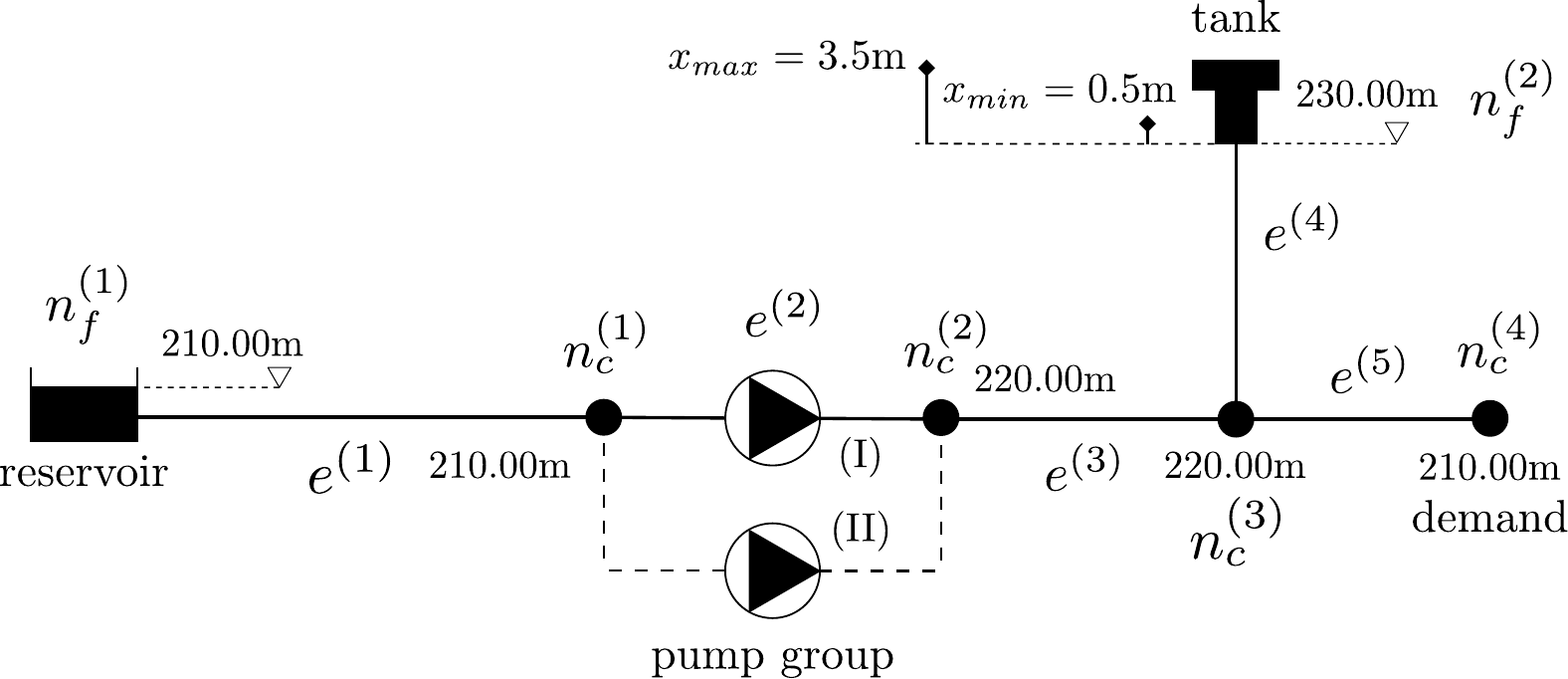}
	\caption{Schematic of a simple network with one tank, one demand point and a single group of two equal \acp{VSP} in parallel.}
	\label{fig:one-tank-network}
\end{figure}

Our method was tested on a model of a simple system illustrated in Fig.~\ref{fig:one-tank-network}.
The network is composed of one fixed-head reservoir, one variable-head tank, two equal parallel \acp{VSP}, 4 pipes and one fixed demand node.
Its purpose is to show the correctness of our method on a simple enough network for which the results are easy to interpret and visualise.

The network was used in two separate analyses.
In the first analysis, the pump schedules were optimized for a 24 hour time horizon for a single default set of inputs and parameters:  tank elevation $z_t= 230$~m, final tank level $x_t^{end}$ = 2.50~m = initial tank level $x_t^{init}$, average demand $\bar{d}$ = 42.7~L/s, tank diameter $D_t$ = 15.00~m.
The aim of the analysis was to study the behaviour of the \ac{MILP} solver on our problem formulation and to verify the correctness of the obtained results.
In the second analysis, a batch of pump schedule optimization was performed for 81 combinations of network parameters and inputs.
The ranges of parameters were as follows: tank elevations $z_t \in [225, 230, 235]$~m, tank level differences $x_t^{end} - x_t^{init} \in [-0.5\,\mathrm{m}, 0.0\,\mathrm{m}, +0.5\,\mathrm{m}]$, average demands $\bar{d} \in [34.16, 42.70, 47.00]$~L/s, tank diameters $D_t\in [12.75, 15.00, 17.25]$ ~m.
The goal was to test the reliability and the robustness of the method under a range of operating points and to measure the calculation times required by CPLEX solver to find optimal solutions.
In both studies, the \ac{MILP} solver was set to terminate upon achieving the MIP gap of 0.05.

\subsection{Results}
\label{subsec:results_case_study_simple}
\begin{figure}[!h]
	\centering
	\includegraphics[width=\linewidth]{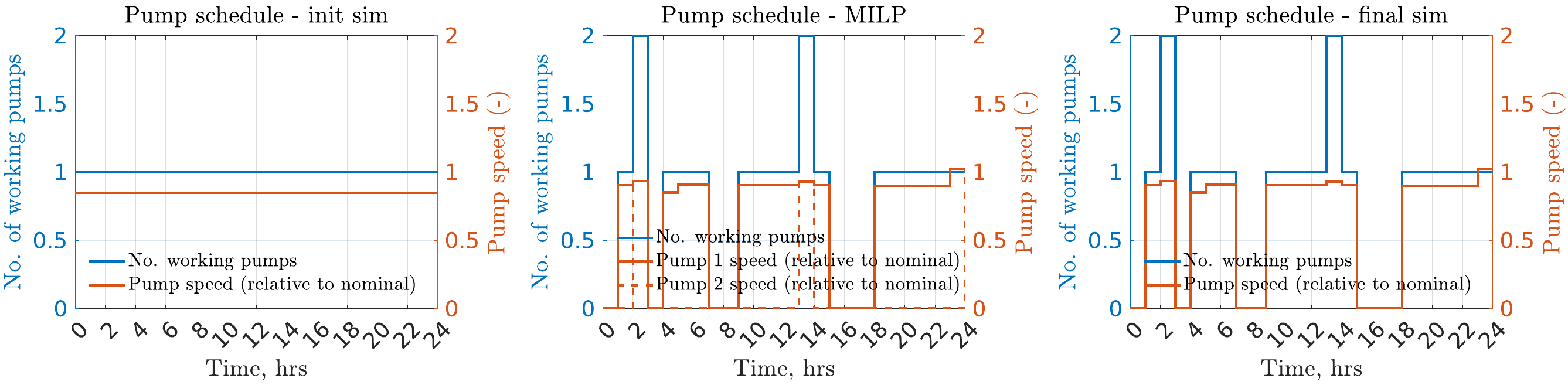}
	\caption{VSP pump schedules.}
	\label{fig:pump_schedules}
\end{figure}

\begin{figure}[!h]
	\centering
	\includegraphics[width=\linewidth]{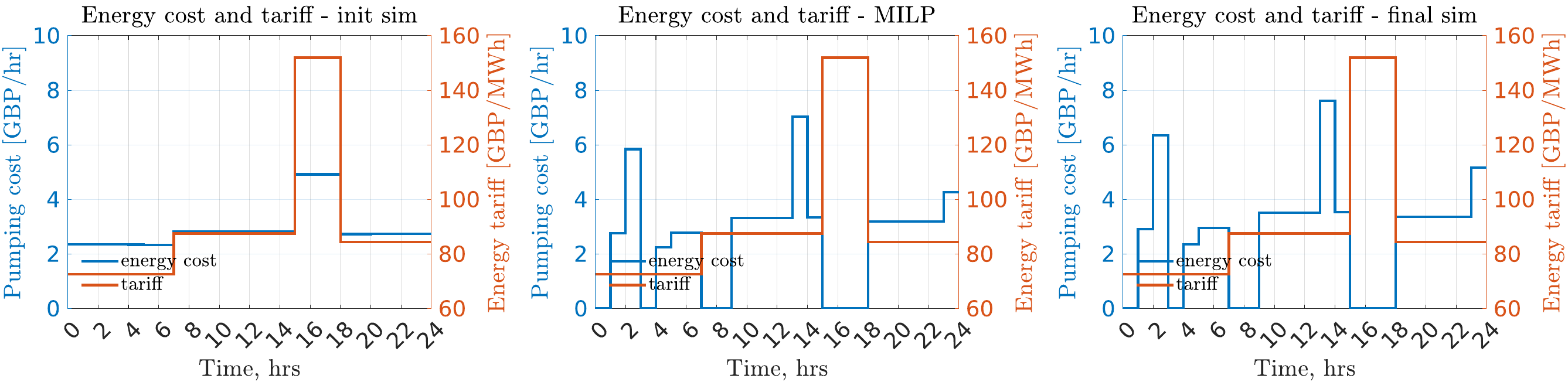}
	\caption{Pumping energy cost and electricity tariff.}
	\label{fig:energy_cost_and_tariff}
\end{figure}

\begin{figure}[!h]
	\centering
	\includegraphics[width=\linewidth]{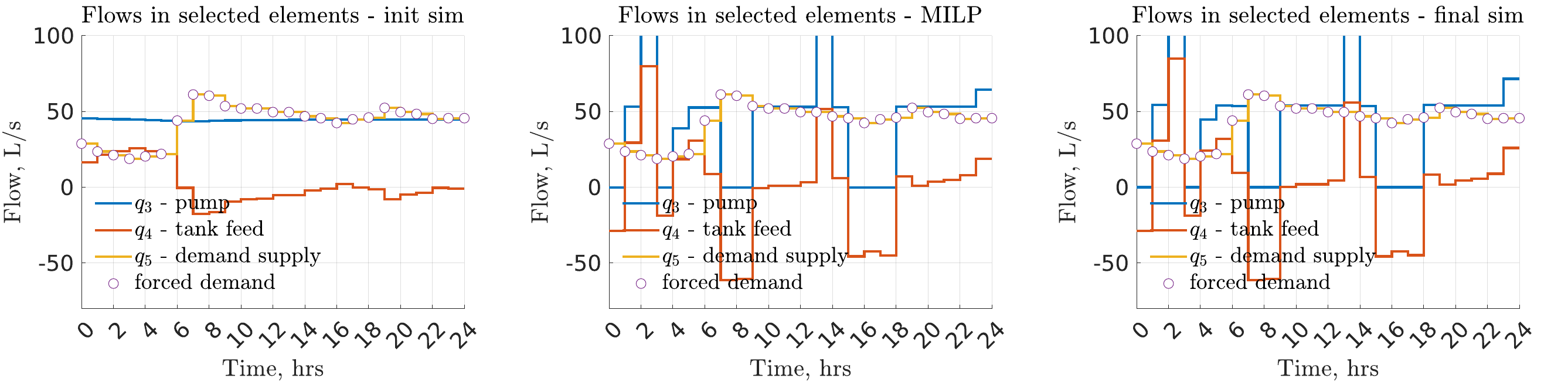}
	\caption{Flows through the pump group, tank feed pipe and demand supply node.}
	\label{fig:flows_selected_elements}
\end{figure}

\begin{figure}[!h]
	\centering
	\includegraphics[width=\linewidth]{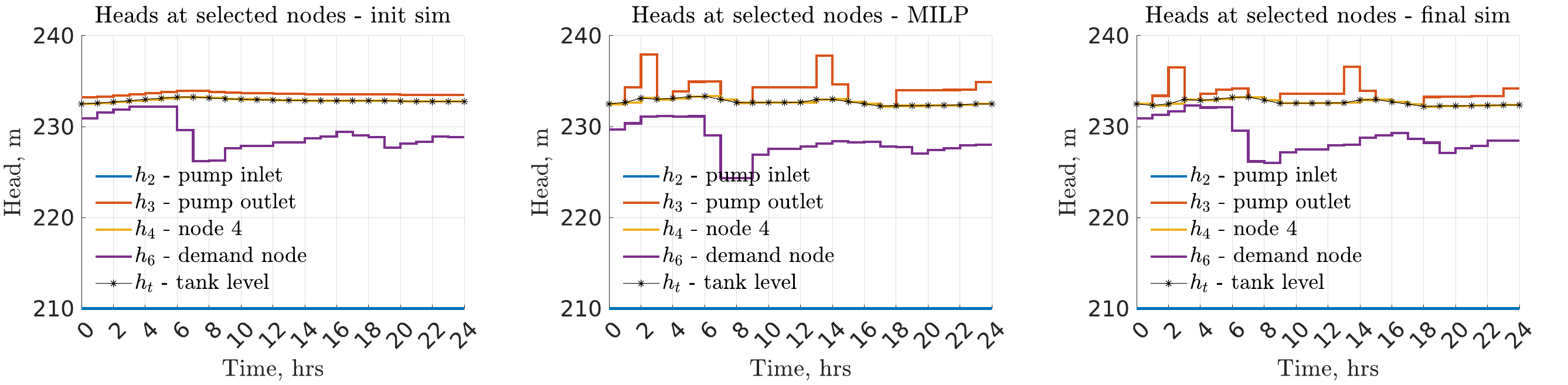}
	\caption{Heads at the pump inlet and outlet, node $n_c^{(3)}$ (node 4), demand node and the tank.}
	\label{fig:heads}
\end{figure}

Results of the initial analysis are shown in Figures~\ref{fig:pump_schedules}, \ref{fig:flows_selected_elements}, \ref{fig:energy_cost_and_tariff} and \ref{fig:heads}.
Results of the initial simulation are shown on the left, the outputs of the \ac{MILP} solver are shown in the middle, whereas the results of the final simulation are shown on the right.
As demonstrated in Fig.~\ref{fig:pump_schedules}, the \ac{MILP} scheduler found an alternative pump schedule to the initial one.
The new schedule reduces the pumps' energy consumption in high tariff periods - see Fig.~\ref{fig:energy_cost_and_tariff}.
Consequently, the operating cost per day of the network reduced from the initial cost of 70.2~GBP to 64.7~GBP (from final simulation), i.e. a 7.8\% reduction.
Fig.~\ref{fig:pump_schedules} (middle) illustrates that switching Pump~1 always precedes switching Pump~2 and the speed of inactive pumps is always set to zero.
It is a desired behaviour enforced by the symmetry breaking constraint and binary linearization of the pump speed, respectively.
The schedule produced by the \ac{MILP} solver is translated into a schedule supported by the simulator, which treats multiple pumps as a group of equal pumps operating at equal speeds, not as separate individual units - see Eqs.~\ref{eq:pump_power_consumption_1} and \ref{eq:pump_head_n_s}.

Flows in the selected network elements and heads in the selected network nodes are significantly altered by the new pump schedule - see the subplots in the right and in the middle vs. the left in Fig.~\ref{fig:flows_selected_elements} and Fig.~\ref{fig:heads}, respectively.
As the pumps' operation switched from constant to one in which the tank's storage capacity is utilized in order to reduce pumping during high tariff periods, the flows in the elements between the pumps and the tank exhibit more variability.
Consequently, the flows are higher, in absolute values, during the times when the tank is filling and emptying.

It is interesting to observe the discrepancies in the network state (heads and flows) between \ac{MILP} that works with linear approximations of network components and the simulator that uses a complete (nonlinear) network model.
We can notice an offset between the head at the demand node $h_6$ returned by the optimizer (middle) and the simulator (right), and relatively higher pump outlet heads $h_3$ returned by the optimizer compared to the simulator.
These differences stem from the inaccuracies introduced by the linear and piece-linear approximations
The approximations can be tightened via introduction of a larger number of piece-linear segments or by iterative adjustment of the locations of the break-points.
Both approaches normally come at the cost of increasing optimization times.
The similarity in tank levels is preserved in \ac{MILP} to a greater degree than the heads in non-storage nodes.
This is due to the fact that tank levels change in response to changes in flows which are affected by approximations less than pressures as the former are forced inputs in demand driven simulations and the latter are the outputs and therefore, are determined by the model.

\begin{figure}
    \centering
    \begin{adjustbox}{valign=t}
    \begin{minipage}{0.49\textwidth}
        \centering
        \includegraphics[width=0.85\textwidth]{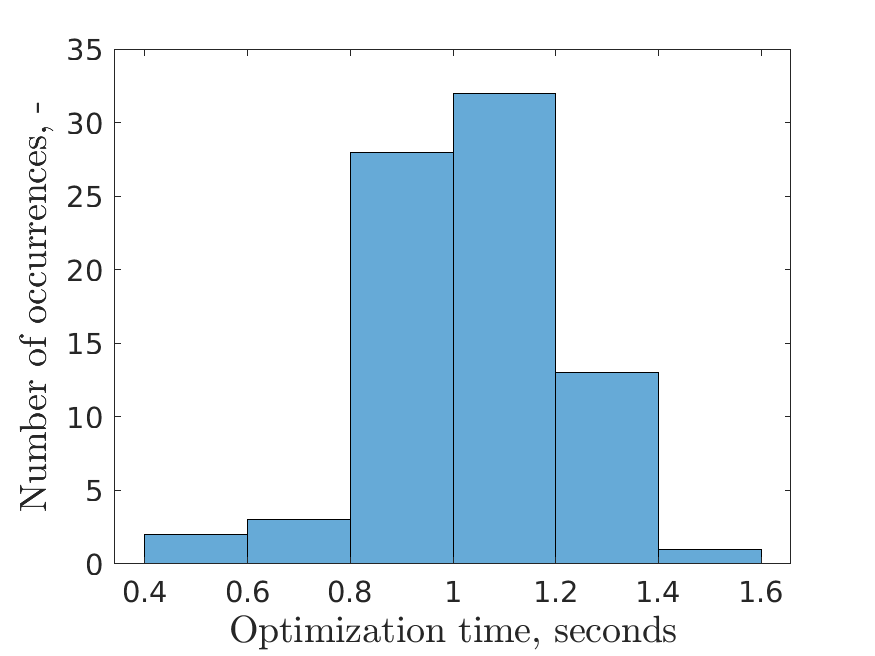} 
        \caption{Histogram of optimization times returned by CPLEX.}
        \label{fig:histogram_running_times}
    \end{minipage}
    \end{adjustbox}\hfill
    \begin{adjustbox}{valign=t}
    \begin{minipage}{0.49\textwidth}
        \centering
        \includegraphics[width=0.85\textwidth]{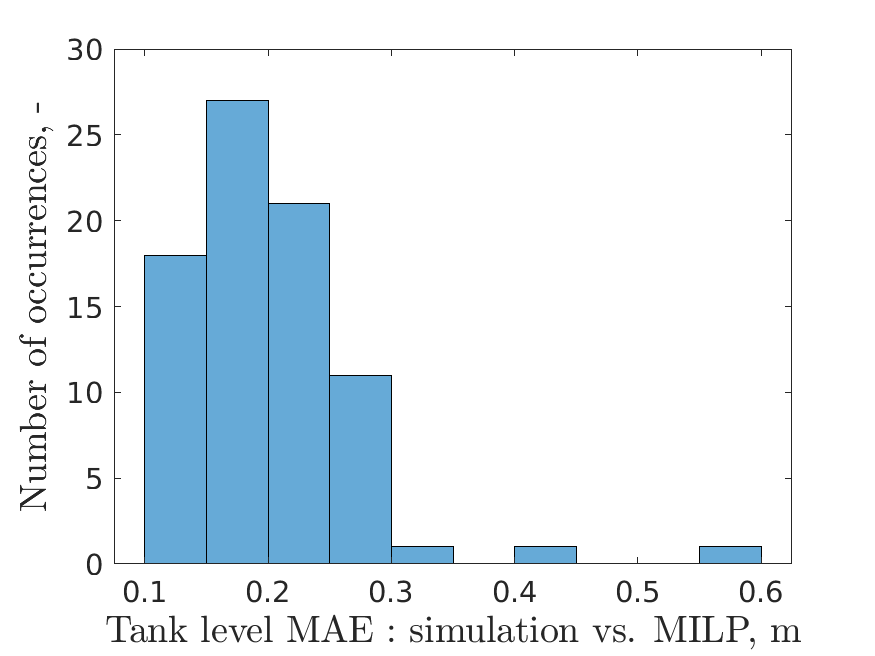} 
        \caption{Histogram of the mean absolute error (MAE) between the simulated tank level and the tank level returned from \ac{MILP}.}
        \label{fig:tank_level_mae}
    \end{minipage}
    \end{adjustbox}
\end{figure}

Fig.~\ref{fig:histogram_running_times} shows the distribution of optimization times out of 80 successful CPLEX runs that produced optimal solutions within MIP gap of 0.05.
All calculations were performed on a laptop equipped with Intel(R) Core(TM) i5-6300U CPU @ 2.40GHz processor and 16GB of RAM and IBM(R) ILOG(R) CPLEX(R) v. 22.1.1.0 \ac{MILP} solver.
The results indicate that the optimization times are bounded (on this problem size) with 75\% of calculations returning optimal solution within $\pm20\%$ of the median time of 1.0 seconds.
Fig.~\ref{fig:tank_level_mae} provides a measure of accuracy of the optimal solution against the simulation.
This metric is quantified as \ac{MAE} between the reservoir level time-series from the simulator and the optimizer.
As demonstrated, 96.6\% of \ac{MILP} outputs return results in which the metric is between 0,1 and 0.3 metres - a rather robust outcome.

\begin{figure}
    \centering
    \begin{minipage}{0.49\textwidth}
        \centering
        \includegraphics[width=1.05\textwidth]{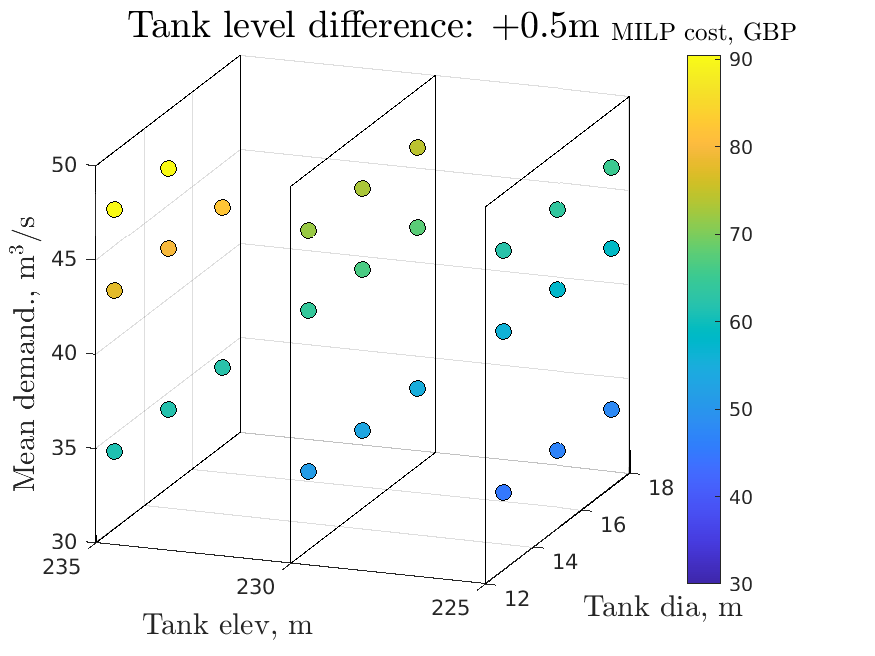} 
        \caption{Optimized pumping costs for a range of tank diameters, tank elevations and water demands, at final tank level difference of +0.5m, i.e. final level 0.5m above the initial level.}
        \label{fig:compound_results_0_5_diff}
    \end{minipage}\hfill
    \begin{minipage}{0.49\textwidth}
        \centering
        \includegraphics[width=1.05\textwidth]{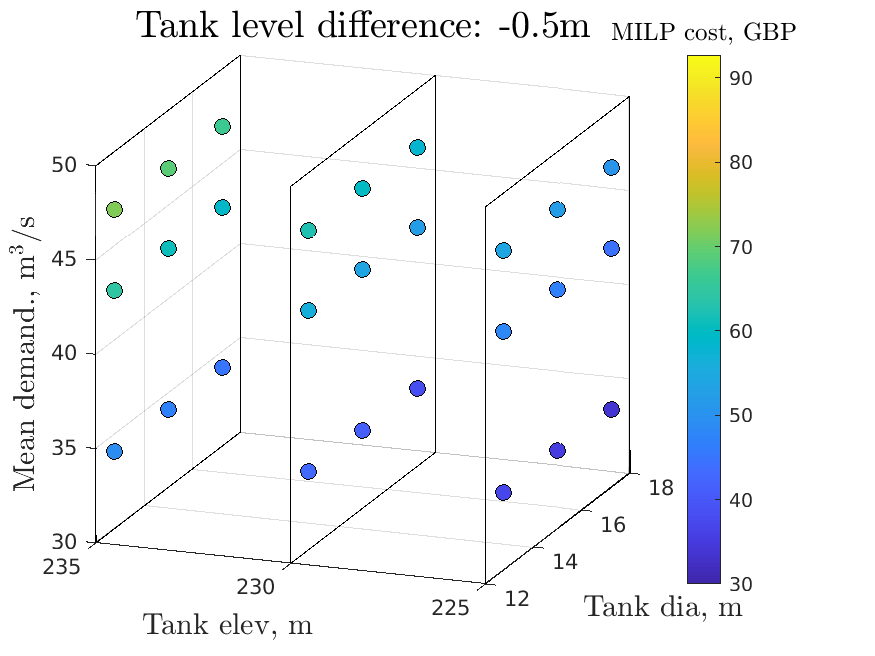} 
        \caption{Optimized pumping costs for a range of tank diameters, tank elevations and water demands, at final tank level difference of -0.5m, i.e. final level 0.5m below the initial level.}
        \label{fig:compound_results_minus_0_5_diff}
    \end{minipage}
\end{figure}

The results of the 53 successful optimizations out of 54 attempts (the other 27 runs for zero tank level difference are not shown) are visualised in Fig.~\ref{fig:compound_results_0_5_diff} and Fig.~\ref{fig:compound_results_minus_0_5_diff}.
The purpose of these visualisations is to demonstrate that the results obtained from the \ac{MILP} solver are physically correct and smooth, which means that the optimizer is able to find a global solution in a robust way for a range of network parameters and inputs.
As suspected, larger operating costs are incurred when the final tank level needs to be 0.5m higher than initial (see Fig.~\ref{fig:compound_results_0_5_diff}), contrary to Fig.~\ref{fig:compound_results_minus_0_5_diff} where the opposite is true.
Operating costs also increase with demand, due higher pumped volumes and larger headlosses, and when the tank is positioned at higher elevations and the pumps need to overcome larger head differences.

\section{Conclusions and further work}
\label{sec:conclusions}
The study demonstrates that mixed integer linear programming with linear and piece-linear approximations of the objective and of the model components is a valid method for finding globally optimal pump schedules in networks containing \acfp{VSP}.
Although the method was tested on a very small network, the average calculation time of approx. 1 second is a promising result indicating that the same approach can be adapted to solving more complex networks.
The method proved to be robust and able to arrive at global optimal solution within similar calculation times for a range of operating points.
As the solvers for mixed integer programming problems, such as CPLEX, GUROBI, MOSEK, etc. have become faster and now support parallel execution, it is perhaps a good idea to start reintroducing mixed integer linear programming techniques, that are known for their stability and robustness, into \acp{WDN} operation studies.
A particularly suited application would be real-time pump optimization.
Another application could be a hierarchical two-level optimization in which the inner pump schedule optimization loop is solved using mixed integer linear programming whilst higher level decisions e.g. long-term policies, design options, etc. are optimized using evolutionary algorithms.
The work presented here shows just one out of many ways of formulating the problem.
It is most likely not the optimal method nor a complete one as it does not support some of the important network elements such as e.g. \acp{PRV}, \acp{CV}, or other aspects such as e.g. handling pressure-dependent demands.
Improvement of solution accuracy and speed can be attempted by experimenting with different component approximation and relaxation techniques, different approximation accuracy improvements and relaxation bound tightening, or different piece-linear approximation representations such as e.g. \acp{SOS}.
These techniques can be borrowed from the existing literature or developed new.
The scalability and the speed of the method can be improved via decomposition techniques such as e.g. Lagrangian or Benders decomposition.
We are currently developing a free open-source Python package that summarizes the state-of-the art in pump scheduling using mixed integer linear programming and can be used for solving practical problems on EPANET networks using different \ac{MILP} formulations and for adding new features and enhancements.
The software (under development) is currently hosted in the \verb|dev-python| branch of the GitHub repository of \emph{MILOPS-WDN - the Mixed Integer Linear Optimal Pump Scheduler} \citep{milops-wdn2023}.

\bibliographystyle{unsrtnat}
\bibliography{references}.

\begin{acronym}
\acro{2D}{two-dimensional}
\acro{ACO}{ant colony optimization}
\acro{ANN}{artificial neural network}
\acro{CV}{check valve}
\acro{EA}{evolutionary algorithm}
\acro{EPS}{extended period simulation}
\acro{FSP}{fixed-speed pump}
\acro{GA}{genetic algorithm}
\acro{GPU}{graphical processing unit}
\acro{HPC}{high performance computing}
\acro{LB}{lower bound}
\acro{LP}{linear programming}
\acro{LPG}{Lagrangian Relaxation and Primal Greedy}
\acro{MILP}{mixed integer linear program}
\acro{MINLP}{mixed integer nonlinear program}
\acro{NRV}{non-return valve}
\acro{MAE}{mean average error}
\acro{MOEA}{multiobjective evolutionary algorithm}
\acro{NLP}{nonlinear programming}
\acro{OR}{operational research}
\acro{OWF}{optimal water flow}
\acro{WDN}{water distribution network}
\acro{PSO}{particle swarm optimization}
\acro{PSP}{pump scheduling problem}
\acro{PRV}{pressure reducing valve}
\acro{SOC}{second order cone}
\acro{SOS}{special ordered set}
\acro{SOS2}{special ordered set type II}
\acro{UB}{upper bound}
\acro{VSP}{variable-speed pump}
\end{acronym}

\section{Appendix}

\subsection{Expanded expressions for $m_j^q$, $m_k^s$ and $c_j$ in the linearized pump power equation \ref{eq:linearized_pump_power_simplified}}
\label{app:lin_pump_power_expressions}

\begin{equation*}
	m_j^q = 3\,a_{3,j}\,q_0^2 + 2\,a_{2,j}\,s_0\,q_0 + a_{1,j}\,s_0^2
\end{equation*}
\begin{equation*}
	m_j^s = a_{2,j}\,q_0^2 + 2\,a_{1,j}\,q_0\,s_0 + 3\,a_{0,j}\,s_0^2
\end{equation*}
\begin{equation*}
	c_j = -2 \left( a_{0,j}\,s_0^3 + a_{1,j}\,q_0\,s_0^2 + a_{2,j}\,q_0^2\,s_0 + a_{3,j}\,q_0^3\right) = - 2\, P_j \, (q_0, s_0)
\end{equation*}

\subsection{Derivation of the linearized pump characteristic}
\label{app:linearized_pump_characteristics_plane_derivation}

The linearized characteristic i.e. the equations of the four planes containing $A_1$, $A_2$, $A_3$ and $A_4$, respectively, are derived from the vector normal to each plane. The normal vector $\vec{n}$ is calculated as a vector product of two vectors formed by the sides of each (triangular) plane.
From definition of orthogonality, the dot product of $\vec{n}$ and $[s,q,H]-p_n$ has to be equal to zero if $[s,n,H]$ lies on the same plane as $p_n$.
Let's follow the procedure for the segment $A_1$.

\begin{equation}
  \vec{n_{A_1}} = \vec{L_1} \times \vec{L_2} =
  \begin{vmatrix}\vec{i} & \vec{j} & \vec{k} \\s_{p_n}-s_{p_1} & q_{p_n}-q_{p_1} & H_{p_n}-H_{p_1} \\
  s_{p_n}-s_{p_2} & q_{p_n}-q_{p_2} & H_{p_n}-H_{p_2}\end{vmatrix}
\end{equation}
\begin{equation}
 \vec{n_{A_1}} \cdot \left( [s,q,H] - [s_n, q_n, H_n] \right) = 0
\end{equation}
Solving the above equation yields the linear plane equation $dd_{j,i} \, s_{j,i} + ee_{j,i} \, q_{j,i} + ff_{j,i} \, H_{j,i} = gg_{j,i}$ for each of the four planes, in which $dd_{j,i}$, $ee_{j,i}$, $ff_{j,i}$ and $gg_{j,i}$ are the calculated coefficients. The procedure is repeated for each of the four segments of the linearized pump characteristic and for each pump $j$.

%


\end{document}